\documentclass[letterpaper, 10pt, conference]{IEEEtran}      

\IEEEoverridecommandlockouts 
\RequirePackage[OT1]{fontenc}
\RequirePackage{amsmath,multicol}
\RequirePackage[colorlinks,citecolor=blue,urlcolor=blue]{hyperref}
\usepackage{algorithm}
\usepackage{algorithmic}
\usepackage{amsthm}
\usepackage[top=1in, left=0.75in, right=0.75in, bottom=0.75in]{geometry}

\usepackage[utf8]{inputenc} 
\usepackage[T1]{fontenc}    
\usepackage{hyperref}       
\usepackage{url}            
\usepackage{booktabs}       
\usepackage{amsfonts}       
\usepackage{nicefrac}       
\usepackage{microtype}      
\usepackage{graphicx}
\usepackage{enumerate}
\usepackage{latexsym}
\usepackage{amssymb}
\usepackage{subfigure}
\usepackage{bm}

\def\bbR{\mathbb{R}}

\def\bbR{\mathbb{R}}
\def\bbP{\mathbb{P}}

\def\kermat{K}

\newcommand{\argmin}{arg\,min}

\newcommand{\la}{\langle}
\newcommand{\ra}{\rangle}

\newcommand{\cH}{\mathcal{H}}

\newcommand{\order}{\ensuremath{\mathcal{O}}}
\newcommand{\Sketch}{\ensuremath{S}}

\newcommand{\numobs}{\ensuremath{n}}

\newcommand{\fhatn}{\widehat{f}_n}

\newcommand{\omegadagger}{\ensuremath{\omega^\dagger}}

\newcommand{\matsnorm}[2]{|\!| #1 | \!|_{{#2}}}

\newcommand{\opnorm}[1]{\ensuremath{\matsnorm{#1}{\tiny{\mbox{op}}}}}
\newcommand{\boldy}{\boldsymbol{y}}

\newtheorem{Theorem}{Theorem}[section]

\newtheorem{Assumption}{Assumption}

\IEEEoverridecommandlockouts                              

\title{
Statistically and Computationally Efficient Variance Estimator for Kernel Ridge Regression\\
}

\author{\IEEEauthorblockN{Meimei Liu}
\IEEEauthorblockA{\textit{Department of Statistical Science} \\
\textit{Duke University}\\
Durham, IN - 27708 \\
Email: meimei.liu@duke.edu}
\and
\IEEEauthorblockN{Jean Honorio}
\IEEEauthorblockA{\textit{Department of Computer Science} \\
\textit{Purdue University}\\
West Lafayette, IN - 47907 \\
Email: jhonorio@purdue.edu}
\and
\IEEEauthorblockN{Guang Cheng}
\IEEEauthorblockA{\textit{Department of Statistics} \\
\textit{Purdue University}\\
West Lafayette, IN - 47907 \\
Email: chengg@purdue.edu}
}

\begin{document}

\maketitle
\thispagestyle{empty}
\pagestyle{empty}

\begin{abstract}
In this paper, we propose a random projection approach to estimate variance in kernel ridge regression.
Our approach leads to a consistent estimator of the true variance, while being computationally more efficient. Our variance estimator is optimal for a large family of kernels, including cubic splines and Gaussian kernels. Simulation analysis is conducted to support our theory.
\end{abstract}

\section{INTRODUCTION}
As a flexible nonparametric tool, kernel ridge regression (KRR) has gained popularity in many application fields, such as machine learning, and visualization; see e.g., \cite{orsenigo2012kernel}. Besides the estimation of the predictive mean, an exploration of the \textit{predictive variance} is also important for statistical inference. Predictive variances can be used for inference, for example, to build confidence intervals; or to select the most informative data points in active learning.
There are two sources of uncertainty in the predictive variance: the noise in the data and the uncertainty in the estimation of the target function. However, calculating the second uncertainty is challenging in KRR on a large data set, since the computational burden increases dramatically with respect to the size of the training set. For example, for $n$ data points, the time and space complexity of kernel ridge regression (KRR) are of $O(n^3)$ and $O(n^2)$ respectively. The above can potentially limit the applicability of KRR to big data scenarios. 

 An efficient way to break the computational bottleneck is low-rank approximation of kernel matrices. Existing methods include dimension reduction (\cite{zhang2005learning,braun2008relevant,balcan2006kernels}), Nystr$\ddot{\textrm{o}}$m approximation 
(\cite{rudi2015less}, \cite{alaoui2015fast}, \cite{musco2016recursive}), and random projections
of large kernel matrices (\cite{yang2015randomized}).
Indeed, low-rank approximation strategies effectively reduce the size of large matrices such that
the reduced matrices can be conveniently stored and processed.

In this paper, we propose a randomly sketched predictive variance to reduce the computational complexity.  Theoretically, we show that given a lower bound of the projection dimension, our approach leads to a consistent estimator of the true variance.
Furthermore, our variance estimator is optimal for a large family of kernel matrices with polynomially and exponentially decaying eigenvalues.
This includes, for instance, cubic splines and Gaussian kernels.

To illustrate the applicability of our theorical contribution, we describe an application of our variance estimator in active learning. In many scenarios, the task of manually labeling (unlabeled) data points is expensive and time-consuming. Therefore it is very important to minimize the number of training examples needed to estimate a particular type of regression function. Suppose we have a set of training examples and labels (responses), and we are permitted to actively choose future unlabeled examples based on the data that we have previously seen. Active learning aims at solving this problem and has been used with various learners such as neural networks ~\cite{mackay1992information, cohn1994neural}, mixture models~\cite{cohn1996active}, support vector machines~\cite{tong2001support} and kernel ridge regression. In active learning, the predictive variance can be viewed as an uncertainty score to iteratively select the most informative unlabeled data points. That is, the largest predictive variance corresponds to the highest uncertainty in $y$ for unlabeled points, which indicates that we may need more information regarding those points.

\section{Preliminaries}

In this section, we introduce kernel ridge regression, its mean prediction and the conditional covariance.
Let $\mathbf{X}=(X_1, \cdots, X_n)^\top $ be the training examples, and $\boldy=(y_1,\cdots, y_n)$ be the corresponding training labels, where $X_i \in \mathcal{X}$ with distribution $P_X$ and $y_i \in \mathcal{R}$ for all $i = 1,\dots,n$. Consider the following nonparametric regression model
\begin{equation}\label{model}
y_i=f^*(X_i) +\epsilon_i, \;\;\; \text{for} \;\;i=1,\cdots, n
\end{equation}
where $\epsilon_i$'s are independent random variables with mean $0$ and variance $\sigma^2$. 
Hereafter, we assume that $\cH $ is a reproducing kernel Hilbert
space (RKHS) associated with a reproducing kernel function $K(\cdot,\cdot)$
defined from $\mathcal{X}\times\mathcal{X}$ to $\bbR$. Let $\langle \cdot,\cdot \rangle_{\cH}$ 
denote the inner product of $\mathcal{H}$ associated with $K(\cdot, \cdot)$, then the reproducing kernel property states that
\begin{equation*}\label{eq:reproducing}
  \la f, K(x,\cdot)\ra_{\cH} = f(x), \;\; \text{for all } f\in \cH.
\end{equation*}
The corresponding norm is defined as $\|f\|_{\cH} : =\sqrt{\la f,f\ra_{\cH}}$ for any $f\in \cH$.

\newgeometry{top=0.75in, bottom=0.75in, left=0.75in, right=0.75in}  
The classic kernel ridge regression (KRR) estimate is obtained via minimizing a penalized likelihood function: 
\begin{equation}
\begin{aligned}\label{eq:krr}
  \fhatn \equiv   \argmin_{f \in \cH} \left\{ \frac{1}{n} \sum_{i=1}^{n} (y_i-f(X_i))^2 +  \lambda  \|f\|_{\cH}^2 \right\}
\end{aligned}
\end{equation}
Let $\kermat$ be the $n$-dimensional kernel matrix with entries $\kermat_{ij}=\frac{1}{n}K(X_i, X_j)$ for $1\le i,j\le n$. By the representer theorem, $\fhatn$  has the form
\begin{equation*}
f(\cdot)=\sum_{i=1}^n \omega_i K(\cdot, X_i)
\end{equation*}
for a real vector $\omega=(\omega_1,\ldots,\omega_n)^\top $, equation (\ref{eq:krr}) reduces to solving the following optimization problem:
\begin{equation}\label{eq:originalkrr}
\omega^\dagger  = \argmin_{\omega \in \bbR^n} \Big \{
\omega^\top  \kermat^2 \omega - \frac{2}{n}\boldy^\top \kermat\omega +  \lambda  \omega^\top  \kermat \omega \Big \}.
\end{equation}
Thus, the KRR estimator is expressed as $\fhatn(\cdot) =  \sum_{i=1}^\numobs \omegadagger_i K(\cdot, X_i)$, where $\omegadagger = \frac{1}{n}(K+  \lambda  I)^{-1}\boldy$. 

For a new testing data point $x$, let $k(x)= (K(x,X_1), K(x,X_2), \cdots, K(x,X_n))^\top $. 
It is easy to calculate its mean prediction  and variance given $\mathbf{X}$ and $\boldy$ as follows:
\begin{align*}
\widehat{y}(x) & = \frac{1}{n}k(x)^\top (K+ \lambda  I )^{-1}\boldy \\
V_1(x) = \text{Var}(\widehat{y}(x)|\mathbf{X},x) &= \frac{\sigma^2}{n^2} k(x)^\top (K+ \lambda  I )^{-2}k(x) .
\end{align*}

Analyzing the conditional variance for the testing data is very important in active learning, since it can act as a guide to select the efficient information we need. However, the time and space taken for solving $(K+ \lambda  I)^{-1}$ is of order $O(n^3)$. This cost is expensive especially when the kernel matrix is dense and the sample size is large.

\section{Randomly Projected Variance}
In this section, we introduce our randomly projected conditional covariance and our main assumptions.
Note that by the Binomial Inverse Theorem, we have 
\begin{equation}\label{eq:binomail}
(K+ \lambda  I )^{-1} = \frac{1}{ \lambda }(I - K( \lambda  K + K^2)^{-1}K).
\end{equation}
To reduce the computational cost, now we propose to replace $(K+ \lambda  I )^{-1}$ by using a randomly projected version as follows
\begin{equation}
\frac{1}{ \lambda }(I -KS^\top ( \lambda  SKS^\top  + SK^2S^\top )^{-1}SK),
\end{equation} 
where $S \in \bbR^{m\times n}$ is a random matrix where each row is independently distributed and sub-Gaussian. Then the conditional variance with the randomly projected matrix can be written as 
\begin{equation}\label{eq:provar}
\begin{aligned}
&V_2(x)  = \text{Var}(\widehat{y}(x)|\mathbf{X},x, S)\\
&=
 \hspace{-2pt}\frac{\sigma^2}{n^2 \lambda ^2} k(x)^\top (I \hspace{-1pt}- \hspace{-1pt}KS^\top ( \lambda  SKS^\top \hspace{-1pt}+\hspace{-1pt}SK^2S^\top )^{-1}SK)^{2}k(x)
\end{aligned}
\end{equation}

The definition of $V_2$ is also our contribution. The variance $V_2$ is different from the variance that could be derived from the results in \cite{yang2015randomized}, which is:  
\begin{equation*}
\begin{aligned}
\hspace{-1pt}V_3(x) = \frac{\sigma^2}{n^2} k(x)^\top K S^\top  ( \lambda  SKS^\top \hspace{-1pt}+\hspace{-1pt}SK^2S^\top )^{-1}SK^2S^\top \\ 
( \lambda  SKS^\top \hspace{-1pt}+\hspace{-1pt}SK^2S^\top )^{-1}SKk(x)
\end{aligned}
\end{equation*}
Unfortunately, understanding the concentration of $V_3$ seems highly nontrivial.
However, our new proposed randomly sketched variance $V_2$ in eq.(\ref{eq:provar}) has nice concentration properties in Theorem \ref{thm:main}.

Note that calculating $( \lambda  SKS^\top  + SK^2S^\top )^{-1}$ only takes the order of $\order(mn^2)$, which enhances the computational efficiency greatly. In Section 4, we will provide a lower bound for $m$ that guarantees stability of the variance after random projection. 

Here we introduce some notations to study the dimension of the random matrix. Define the efficiency dimension as 
\begin{align}\label{eq:testing dimen}
s_ \lambda  = \argmin \{j:\widehat{\mu}_j \leq  \lambda\}-1,
\end{align}
where $\widehat{\mu}_j$ is the $j$-th highest eigenvalue of the kernel matrix $\kermat$. More formally, let $\kermat=UDU^\top $, where $U \in \bbR^{n\times n}$ is an orthonormal matrix, i.e., $UU^\top $ is an $n\times n$ identity matrix, and $D\in \bbR ^{n\times n}$ is a diagonal matrix  with diagonal elements $\widehat{\mu}_1 \geq \widehat{\mu}_2 \geq \cdots \geq \widehat{\mu}_n >0 $. 

In this paper, we consider random matrices with independent sub-Gaussian rows. For the random matrix $\Sketch$, the $i^{th}$ row $S_i\in \bbR^n$ is sub-Gaussian if for all $u\in \bbR^n$, $\langle S_i, u\rangle$ are sub-Gaussian random variables, i.e.,
$$
\bbP\{|\langle S_i, u\rangle|>t\} \leq e\cdot \exp\{-t^2\}.
$$
Matrices fulfilling the above condition include all matrices with independent sub-Gaussian entries as a particular instance. The class of sub-Gaussian variates includes for instance Gaussian variables, any bounded random variable (e.g. Bernoulli, multinomial, uniform), any random variable with strongly log-concave density, and any finite mixture of sub-Gaussian variables. In the following of the paper,  we scale the random matrix by $\sqrt{m}$ for analyzing convenience. 

Next, we state our main assumption and some useful results related to the randomly projected kernel matrix.

\begin{Assumption}\label{asmp:random1}
Let $S$ be a sub-Gaussian random matrix with independent rows. Let $ \lambda  \rightarrow 0$ and $ \lambda  \gg 1/n$. Set the projection dimension $m\geq d s_ \lambda $, where $d$ is an absolute constant. For $K=UDU^\top $, let $U=(U_1, U_2)$ with $U_1 \in \bbR^{n\times s_ \lambda }$, and $U_2 \in \bbR^{n \times (n-s_ \lambda )}$. Let $D=\begin{pmatrix}
  D_1 & 0  \\
  0 & D_2 
 \end{pmatrix}
 $ with $D_1 \in \bbR^{s_ \lambda  \times s_ \lambda }$
and $D_2 \in \bbR^{(n-s_ \lambda ) \times (n-s_ \lambda )}$.
We assume that $S$ satisfies the following conditions: 
\begin{enumerate}[(i)]
\item $1/2 \leq  \lambda _{\min}(SU_1) \leq  \lambda _{\max}(SU_1) \leq 3/2$ with probability greater than $1-2\exp\{-cm\}$, where $c$ is an absolute constant independent of $n$. 
\item $\opnorm{SU_2D_2^{1/2}} \leq c' \lambda ^{1/2}$ with probability greater than $1-2\exp\{-c^{''}m\}$, where $c'$ and $c^{''}$ are constants independent of $n$. 
\end{enumerate}
\end{Assumption}
In Assumption A1, the kernel matrix is partitioned into a summation of the form $K= U_1 D_1 U_1^\top  + U_2 D_2 U_2^\top $ where $U_1$ contains the first $s_ \lambda $ columns of the orthonormal matrix $U$, which correspond to the first leading eigenvalues of the kernel matrix; and $U_2$ contains the rest of the $n-s_ \lambda $ columns of the $U$, which correspond to the smallest $n-s_ \lambda $ eigenvalues. In most cases, the smallest $n-s_ \lambda $ eigenvalues are neglectable due to a fast decaying rate of the eigenvalues. Assumption \ref{asmp:random1} (i) ensures that the randomly projected eigenvectors corresponding to the leading eigenvalues still preserve the distance between each other approximately; Assumption \ref{asmp:random1} (ii) ensures that the operator norm of the lowest ``neglectable'' part would not change too much after the random projection. Random matrices satisfying Assumption A1 include sub-Gaussian random matrices (see detailed proof in \cite{liu2018nonparametric}), as well as matrices constructed by randomly sub-sampling and rescaling the rows of a fixed orthonormal matrix. We refer the interested reader to \cite{yang2015randomized}, \cite{vershynin2010introduction} for more details.

Our work differs from \cite{yang2015randomized} in several fundamental ways. \cite{yang2015randomized} focuses on the (mean) prediction error on a training set, and it is unclear how this relates to a prediction error on a testing set. In contrast to \cite{yang2015randomized}, we target the variance of the prediction error, and focus on prediction on a test set. Additionally, note that we define $s_\lambda$ by the tuning parameter $\lambda$ as in eq. (\ref{eq:testing dimen}), which is different from \cite{yang2015randomized}.

\section{Main Results}
Recall that for a new testing data $x$, the conditional variance given training data $\mathbf{X}$ and $\boldy$ is $V_1(x)$. 
In this section, we will show that for the new data $x$, our new proposed randomly sketched conditional variance $V_2(x)$ can provide a stable approximation for the original conditional covariance $V_1(x)$.
\begin{Theorem}\label{thm:main}
Under Assumption \ref{asmp:random1}, suppose $ \lambda  \rightarrow 0$ as $n\to \infty$, $ \lambda  \gg n^{-1}$, and the projection dimension $m \geq d s_ \lambda $. Then with probability at least $1-2\exp(-cm)$, with respect to the random choice of $S$, we have
$$
\sup_{x\in\mathcal{X}}|V_1(x)-V_2(x)| \leq \frac{c'\sigma^2}{n \lambda },
$$
where $c$ and $c'$ are absolute constants independent of $n$. 

\end{Theorem}

As shown in Theorem \ref{thm:main}, the convergence rate involves $\lambda$ directly. Normally, we choose $\lambda$ as the optimal one to achieve minimax optimal estimation. 
Next, we provide some examples to show how to choose the lower bound of the projection dimension for the random matrix $S$ and the corresponding optimal $\lambda$.

\begin{proof}
 Let $k(x)=(g(X_1), \cdots, g(X_n))^\top $ with $g(X_i) = K(x, X_i)= \la K(x,\cdot), K(\cdot,X_i) \ra$, then we have that 
$g(\cdot) = K(x,\cdot) \in \cH$. We denote $g^*= k(x)=(g(X_1), \cdots, g(X_n))^\top $, and let $\sigma=1$. Then $V_1(x)$ can be written as 
\begin{equation*}
V_1(x) = \frac{\sigma^2}{n^2 \lambda ^2}\|(I-K( \lambda  K + K^2)^{-1}K)g^*\|_2^2,
\end{equation*}
by  eq.(\ref{eq:binomail}), where $\|\cdot\|_2$ is the Euclidean norm. Furthermore
\begin{equation*}
V_2(x) = \frac{\sigma^2}{n^2 \lambda ^2}\|(I-KS^\top ( \lambda SKS^\top +SK^2S^\top )^{-1}SK)g^*\|_2^2.
\end{equation*}
and therefore:
\begin{align}
 & \sup_{x\in \mathcal{X}}|V_1(x)-V_2(x)| \nonumber\\
 = &\sup_{x\in \mathcal{X}} \frac{\sigma^2}{n^2 \lambda ^2}[(I-K( \lambda  K + K^2)^{-1}K)g^* \nonumber\\
 & + (I-KS^\top ( \lambda SKS^\top +SK^2S^\top )^{-1}SK) g^* ]^\top   \nonumber\\
& \cdot [(I-K( \lambda  K + K^2)^{-1}K)g^* \nonumber\\
& - (I-KS^\top ( \lambda SKS^\top +SK^2S^\top )^{-1}SK) g^* ] \nonumber\\
\leq &  \frac{\sigma^2}{n \lambda ^2} (T_1 +T_2)^2
\end{align}\label{diff:V1V2}
where in the last step, we used Cauchy Schwarz inequality and the triangle inequality. In the above, $T_1$ and $T_2$ are defined as follows:
\begin{align*}
T_1 & = \frac{1}{\sqrt{n}} \|g^* - K( \lambda  K + K^2)^{-1}K g^*\|_2\\
T_2 & = \frac{1}{\sqrt{n}} \|g^* - KS^\top ( \lambda SKS^\top +SK^2S^\top )^{-1}SK g^*\|_2
\end{align*}

Next, we prove that 
\begin{equation}\label{ineq: T1T2bound}
T_1^2 \lesssim  \lambda ,\quad   \quad T_2^2 \lesssim  \lambda .
\end{equation} 

Before the proof of eq.(\ref{ineq: T1T2bound}), we first consider the optimization problem 
\begin{equation}\label{opt:random}
\widehat{\alpha}=\argmin_{\alpha \in \bbR^{m}} \frac{1}{n}\|g^* -nKS^\top  \alpha\|_2^2 +  n\lambda  \|K^{1/2}S^\top \alpha\|_2^2,
\end{equation}
which has the solution $\widehat{\alpha} =\frac{1}{n}( \lambda SKS^\top +SK^2S^\top )^{-1}SK g^*$. In this case $T_2^2 = \frac{1}{n}\|g^* -nKS^\top  \widehat{\alpha}\|_2^2$.

Therefore, to prove $T_2^2 \leq  \lambda $, we only need to find a vector $\widetilde{\alpha}$, such that 
$\frac{1}{n} \|g^* -n KS^\top  \widetilde{\alpha}\|_2^2 +  n\lambda  \|K^{1/2}S^\top \widetilde{\alpha}\|_2^2 \leq c_1  \lambda $. This will imply
\begin{align}\label{ineq:construct}
& \frac{1}{n}\|g^* -nKS^\top  \widehat{\alpha}\|_2^2 + n \lambda  \|K^{1/2}S^\top \widehat{\alpha}\|_2^2 \nonumber\\
\leq & \frac{1}{n}\|g^* -nKS^\top  \widetilde{\alpha}\|_2^2 +  n\lambda  \|K^{1/2}S^\top \widetilde{\alpha}\|_2^2 \leq c_1  \lambda . 
\end{align}

By definition of $s_ \lambda $, when $1\leq j \leq s_ \lambda $ then $\widehat{\mu}_j \geq  \lambda $ and when $s_ \lambda < j \leq n$ then $\widehat{\mu}_j \leq  \lambda $. Let $g^*=(g^*_1, g^*_2)$, where $g^*_1 \in \bbR^{s_ \lambda }$, and $g^*_2 \in \bbR^{n-s_ \lambda }$. Let $z=\frac{1}{\sqrt{n}}U^\top g^* = (z_1,z_2)$ correspondingly.  Also, divide $D$ into $D_1, D_2$, where $D_1$, $D_2$ are $s_ \lambda  \times s_ \lambda $  and $(n-s_ \lambda )\times (n-s_ \lambda )$ dimension diagonal matrix, respectively. Let $\widetilde{S} = (\widetilde{S}_1, \widetilde{S}_2)$, with $\widetilde{S}_1\in \bbR^{s\times s_ \lambda }$ as the left block and $\widetilde{S}_2\in \bbR^{s\times (n-s_ \lambda )}$ as the right block.  
We construct a vector $\widetilde{\alpha}$ by setting $\widetilde{\alpha} = \frac{1}{\sqrt{n}}\widetilde{S}_1 (\widetilde{S}_1^\top  \widetilde{S}_1)^{-1}D_1^{-1}z_1 \in \bbR^{s}$.  By plugging $\widetilde{\alpha}$ into eq.(\ref{opt:random}), we have that 
\begin{align*}
& \frac{1}{n}\|g^*-nUD\widetilde{S}^\top  \widetilde{\alpha}\|_2^2  \\
= & \|z_1 - \sqrt{n}D_1 \widetilde{S}_1^\top  \widetilde{\alpha}\|_2^2 + \|z_2 - D_2 \widetilde{S}_2^\top \widetilde{S}_1(\widetilde{S}_1^\top \widetilde{S}_1)^{-1}D_1^{-1}z_1\|_2^2\\
 = &  G_1^2 + G_2^2.
\end{align*}
Clearly, in our construction $G_1^2 =0$, and thus we focus on analyzing $G_2$. For any $g(\cdot) \in \cH$, there exists a vector $\beta \in \bbR^n$, such that $g(\cdot) = \sum_{i=1}^n K(\cdot,X_i)\beta_i + \xi(\cdot)$, where $\xi(\cdot)\in \cH$, and such that $\xi$ is orthogonal to the span of $\{K(\cdot, X_i), i=1,\cdots,n\}$. Therefore, $\xi(X_j)=\la \xi, K(\cdot,X_j) \ra =0$, and $g(X_j)=\sum_{i=1}^n K(X_i,X_j)\beta_i$. Thus $g^* = nK\beta$, where $K$ is the empirical kernel matrix. Assume that $\|g\|_{\cH} \leq 1$, then 
\begin{align*}
& n\beta^\top  K \beta  \leq 1 \quad \Rightarrow n\beta^\top K K^{-1} K\beta^\top   \leq 1 \quad \\
&\Rightarrow \frac{1}{n}g^* K^{-1}g^* \leq 1 \quad  \Rightarrow  \frac{1}{n}g^*UD^{-1}U^\top  g^* \leq 1
\end{align*}
Then, we have the ellipse constraint that $\|D^{-1/2}z\|_2 \leq 1$, where $z=\frac{1}{\sqrt{n}}U^\top g^*$. 

Since we have $\|D_1^{-1/2}z_1\|_2 \leq 1$, $\|D_2^{-1/2}z_2\|_2 \leq 1$, which implies $g^{*T}U_2U_2^\top  g^* \leq n \lambda $, we have that

\begin{align*}
G_2  \leq & \|z_2\|_2 + \opnorm{\sqrt{D_2}}\opnorm{\sqrt{D_2}\widetilde{S}_2^\top }\opnorm{\widetilde{S}_1} \opnorm{(\widetilde{S}_1^\top \widetilde{S}_1)^{-1}}\\
&\cdot \opnorm{D_1^{-1/2}}\opnorm{D_1^{-1/2}z_1} \leq c\sqrt{ \lambda }
\end{align*}
Therefore, we have $\|z-\sqrt{n}D\widetilde{S}^\top  \widetilde{\alpha}\|_2^2 \leq c' \lambda $. For the penalty term, 
\begin{align*}
& n\widetilde{\alpha}^\top SKS^\top  \widetilde{\alpha}  \leq z_1^\top  D_1^{-1}z_1 + \|z_1^\top D_1^{-\frac{1}{2}}\|_2\opnorm{D_1^{-\frac{1}{2}}} \\
 & \cdot \|\widetilde{S}_2\sqrt{D_2}\|\opnorm{\sqrt{D_2}\widetilde{S}^\top } \opnorm{D_1^{-\frac{1}{2}}}\opnorm{D_1^{-\frac{1}{2}}z_1} \leq c'',
\end{align*}
where $c''$ is a constant. Finally, by eq.(\ref{ineq:construct}), we can claim that 
$$
\frac{1}{n}\|KS^\top ( \lambda SKS^\top +SK^2S^\top )^{-1}SK g^* -g^*\|_2^2 \leq c_1 \lambda  ,
$$
where $c_1$ is some constant.

Similarly, to prove $T_1^2\lesssim  \lambda $, we can treat $S$ as an identity matrix. Consider the following optimization problem 
\begin{equation*}\label{opt:orginal}
\widehat{w}=\argmin_{w \in \bbR^{n}} 
\frac{1}{n}\|g^* -nK w\|_2^2 + n \lambda  \|K^{1/2}w\|_2^2,
\end{equation*}

which has the solution $\widehat{w} =\frac{1}{n}( \lambda K+K^2)^{-1}K g^*$. In this case $T_1^2 = \frac{1}{n}\|g^*-nK\widehat{w}\|_2^2$. Therefore, we only need to find a $\widetilde{w}$, such that 
$ \frac{1}{n}\|g^* -nK \widetilde{w}\|_2^2 +  n\lambda  \|K^{1/2}\widetilde{w}\|_2^2 \leq c_2  \lambda $. This will imply 
\begin{equation}\label{ineq:2construct}
\begin{aligned}
 & \frac{1}{n}\|g^* -nK \widehat{w}\|_2^2 + n \lambda  \|K^{1/2}\widehat{w}\|_2^2 \\
\leq & \frac{1}{n}\|g^* -nK\widetilde{w}\|_2^2 + n \lambda  \|K^{1/2}\widetilde{w}\|_2^2 \leq c_2  \lambda .
\end{aligned}
\end{equation}
Here we construct $\widetilde{w} = \frac{1}{\sqrt{n}}U_1D_1^{-1}z_1$, 
\begin{align*}
& \frac{1}{n}\|g^*-UDU^\top  \widetilde{w}\|_2^2 \\
 =& \|z_1 - D_1 U_1^\top  \widetilde{w}\|_2^2 
+ \|z_2 - D_2 U_2^\top U_1(U_1^\top U_1)^{-1}D_1^{-1}z_1\|_2^2\\
 = & \|z_2\|_2^2 \leq c_2 \lambda .
\end{align*}
For the penalty term, $n\widetilde{w}^\top K \widetilde{w} = z_1^\top D_1^{-1}z_1 \leq 1$. Therefore, combining with eq.(\ref{ineq:2construct}), we have 
$$
\frac{1}{n}\|K( \lambda  K + K^2)^{-1}K g^* -g^*\|_2^2 \leq c_2  \lambda .
$$

Finally, by eq.(4.1) , we have 
$$
\sup_{x\in\mathcal{X}}|V_1(x) - V_2(x)| \leq \frac{\sigma^2}{n  \lambda ^2} (T_1+ T_2)^2 \lesssim \frac{\sigma^2}{n \lambda }
$$

 \end{proof}

\textbf{Example 1}: Consider kernels with polynomially decaying eigenvalues $\mu_k \asymp k^{-2\alpha}$ for $\alpha\geq 1$. Such kernels include the $\alpha-$order periodic Sobolev space, for $\alpha=2$, which corresponds to the cubic spline. Since the optimal rate of $ \lambda $ to achieve the minimax estimation error is of order $n^{-\frac{2\alpha}{2\alpha+1}}$, we get the corresponding optimal lower bound for the projection dimension $m \gtrsim s_ \lambda \asymp n^{\frac{1}{2\alpha+1}}$. Furthermore, the difference between original conditional variance and randomly sketched conditional variance $|V_1(x)-V_2(x)|$ can be bounded by the order of $O(n^{-\frac{1}{2\alpha+1}})$. 

\textbf{Example 2}: Consider kernels with exponentially decaying eigenvalues $\mu_k \asymp e^{-\alpha k^p}$ for $p>0$, which include the Gaussian kernel with $p=2$. Since the optimal rate of $ \lambda $ to achieve the minimax estimation rate is of order $(\log n)^{1/p}/n$, we get the corresponding lower bound $m\geq s_ \lambda  \asymp (\log (n))^{1/p}$, and $|V_1(x)-V_2(x)| \lesssim (\log n)^{-1/p}$.

\section{Experiments}

In this section, we verify the validity of our theoretical contribution (Theorem \ref{thm:main}) through synthetic experiments.

Data were generated based on eq.(\ref{model}) with the predictor $X$ following a uniform distribution on $[0,1]$, $f^{\ast}(x) = -1 + 2 x^2$, and $\epsilon_i\sim N(0,\sigma^2)$. We used Gaussian random projection matrices.

For the polynomial kernel, Figure 1 (a) shows the gap $\sup_{x\in[0,1]} |V_1(x)-V_2 (x)|$ with the training sample size $n$ ranging from $50$ to $1000$, while fixing $\sigma = 1$, and the projection dimension $m=\lceil 1.5n^{1/(2\alpha+1)}\rceil$ with $\alpha=2$. Note that with the increase of the sample size $n$, the gap $|V_1-V_2|$ decreases with the rate $O(1/n)$ as predicted by Theorem \ref{thm:main}. For Figure 1 (c), we fix the sample size as $n=1000$, and $\sigma =1$, while varying the projection dimension $m=\lceil 1.2n^{c/(2\alpha+1)}\rceil$ with $c$ ranging from $0.4$ to $1.9$. Note that an increase of $m$ leads to a smaller gap $\sup_{x} |V_1(x)-V_2 (x)|$, but this improvement is no longer obvious when $m\geq 1.2n^{1/(2\alpha+1)}$ (or equivalently when c = 1), which is the optimal projection dimension demonstrated in Example 1. In Figure 1 (e), we vary $\sigma$ from $0.5$ to $5$, while fixing sample size $n=1000$ and projection dimension as $\lceil 1.5n^{1/(2\alpha+1)}\rceil$ with $\alpha=2$; Note that $\sup_{x} |V_1(x)-V_2 (x)|$ increases almost linearly with respect to $\sigma$, which is consistent with our theory. 

For the Gaussian kernel, Figure 1 (b) shows gap $\sup_{x\in[0,1]} |V_1(x)-V_2 (x)|$ with the training sample size $n$ ranging from $50$ to $1000$, while fixing $\sigma = 1$ and the projection dimension $m = \lceil 2\sqrt{\log(n)}\rceil$. Note that the gap $|V_1-V_2|$ decreases with the rate $O(1/n)$ as predicted by Theorem \ref{thm:main}. For Figure 1 (d), we fix the sample size as $n=1000$, and $\sigma = 1$, while varying the projection dimension $m=\lceil 1.2(\log(n))^{c/2}\rceil$ with $c$ ranging from 0.3 to 1.8. Note that an increase of $m$ leads to a smaller gap $\sup_{x} |V_1(x)-V_2 (x)|$, but this improvement is no longer obvious when $m\geq 1.2 (\log (n))^{1/2}$ (or equivalently when c = 1), which is the optimal projection dimension demonstrated in Example 2. In Figure 1 (f), we fix $n=1000$ and $m= \lceil 2\sqrt{\log(n)}\rceil$, but vary $\sigma$ from $0.5$ to $5$. As in the previous experiment, note that $\sup_{x}|V_1(x)-V_2(x)|$  increases almost linearly with respect to $\sigma$, which is consistent with our theory.

In Appendix \ref{sec:append:sim}, we show additional synthetic experiments verifying our theoretical contribution. We further illustrate the use of our projected variance estimator in active learning, in synthetic data as well as two real-world datasets. In this illustrative application, by using the randomly sketched predictive variance, the computational complexity is reduced from $O(n^3)$ to $O(m n)$, where $m$ is the projection dimension. Given our theoretical finding, our variance estimator does not sacrifice statistical accuracy.

\begin{figure}[thpb]\label{fig:compare:both}
  \centering
     \begin{tabular}{cc}
      \includegraphics[scale=0.32]{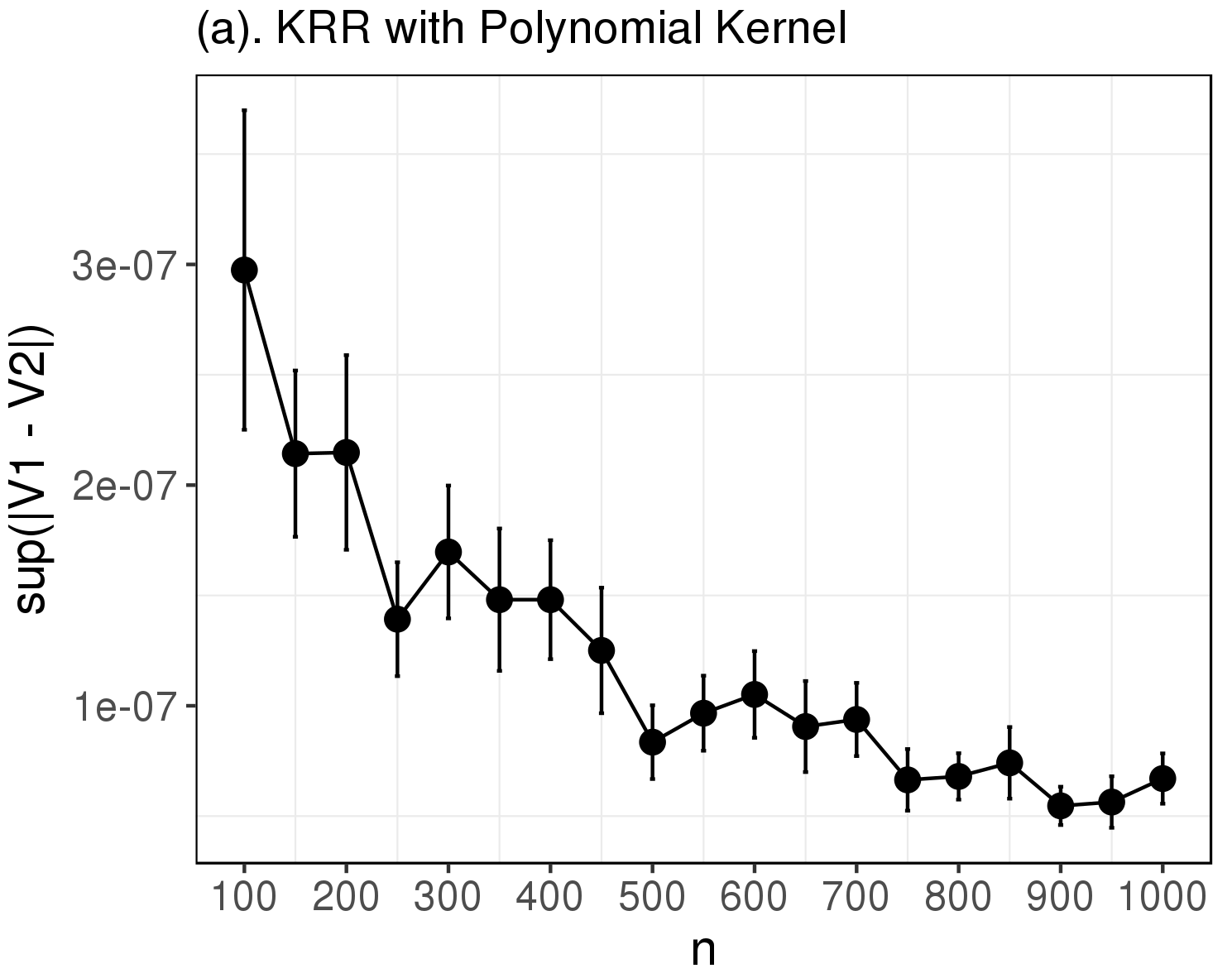}&\includegraphics[scale=0.32]{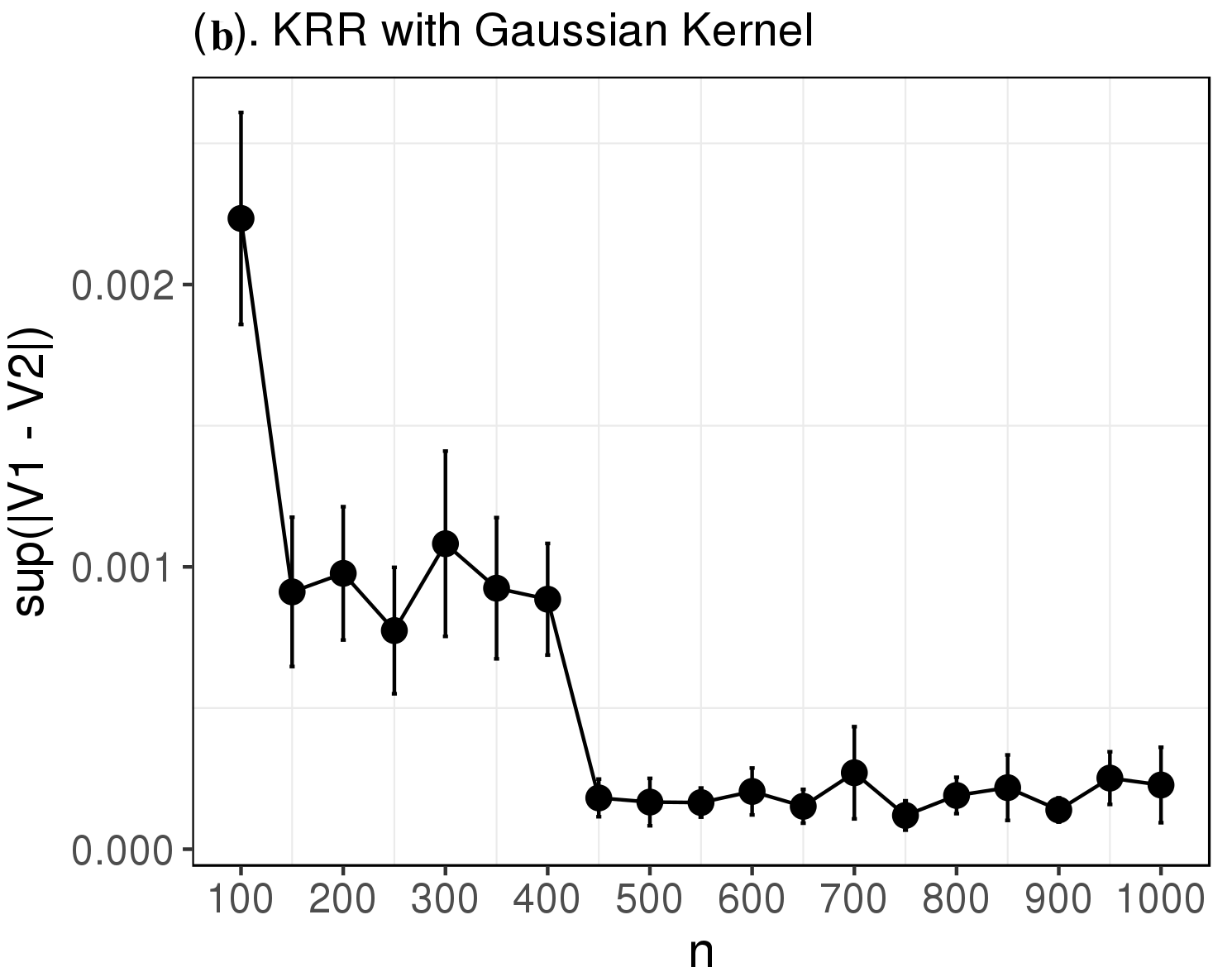}\\
    \includegraphics[scale=0.32]{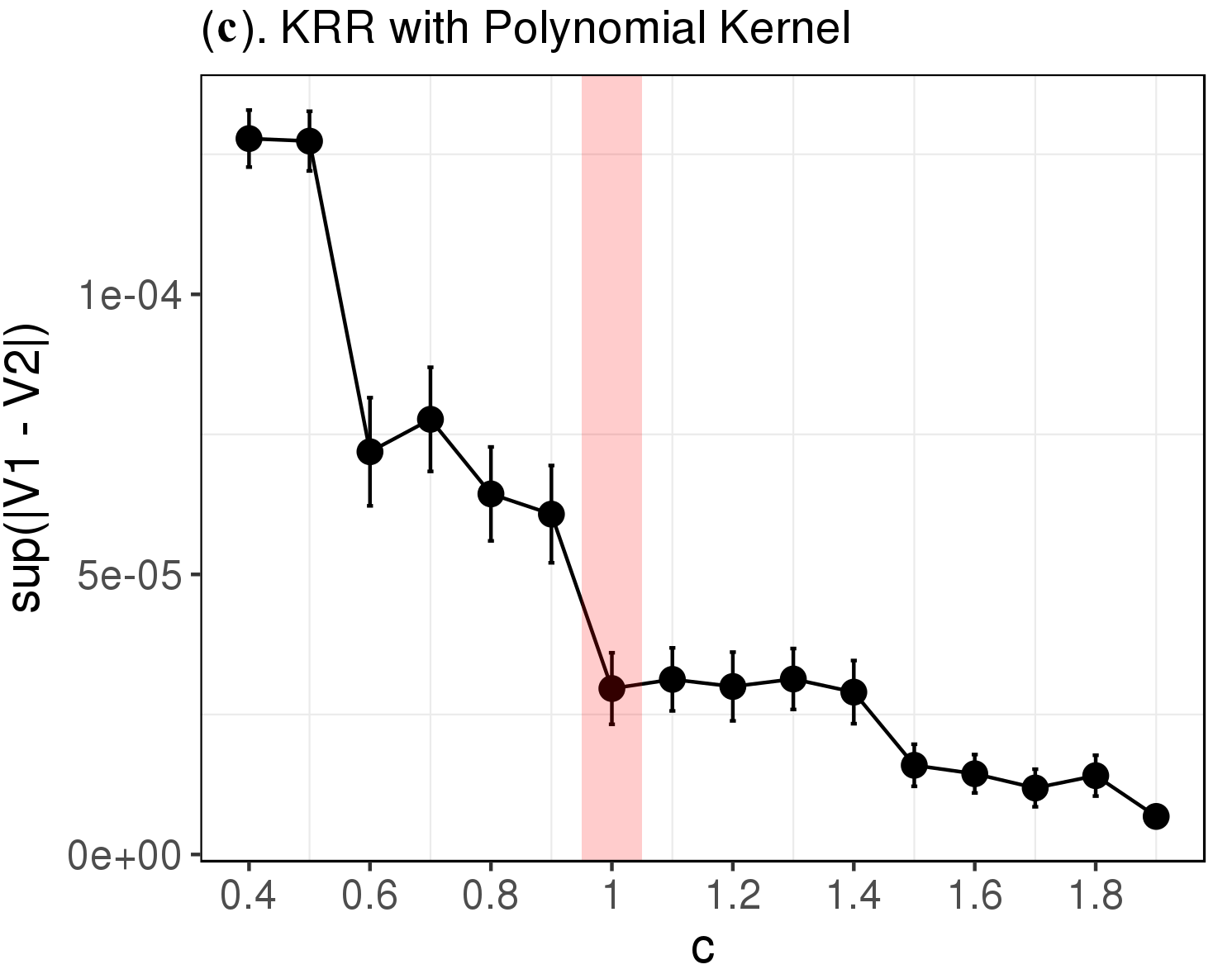}&\includegraphics[scale=0.32]{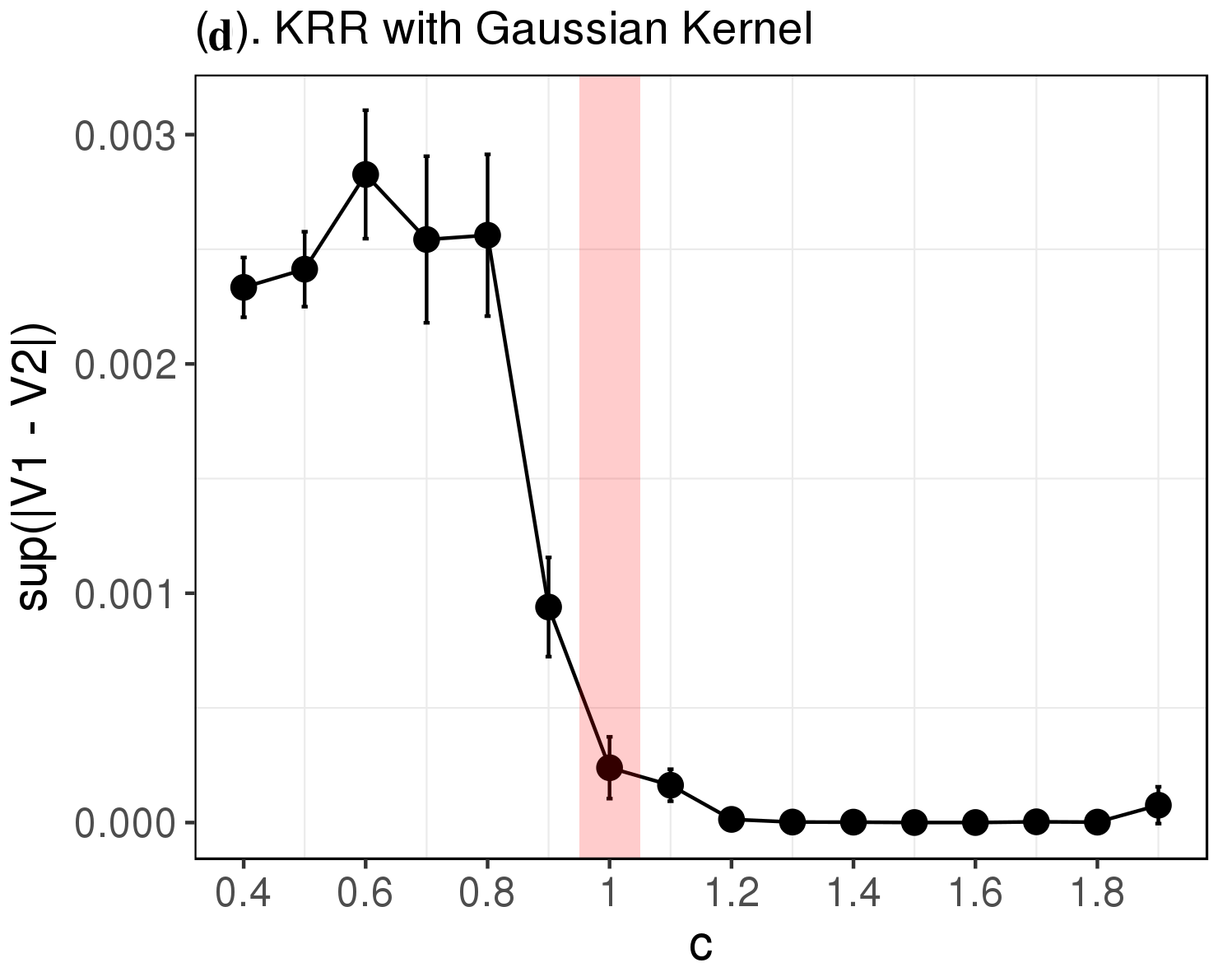}\\
    \includegraphics[scale=0.32]{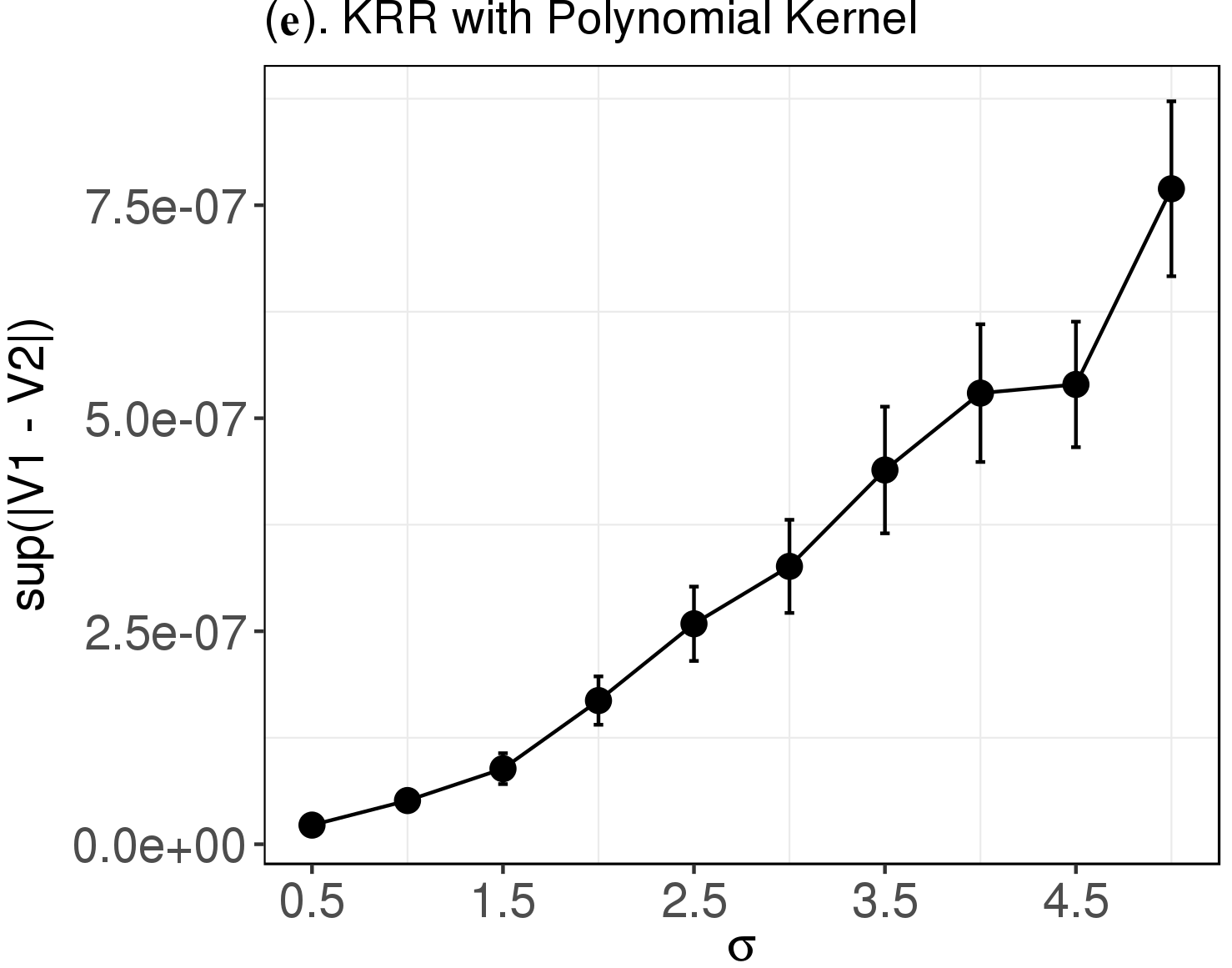}&\includegraphics[scale=0.32]{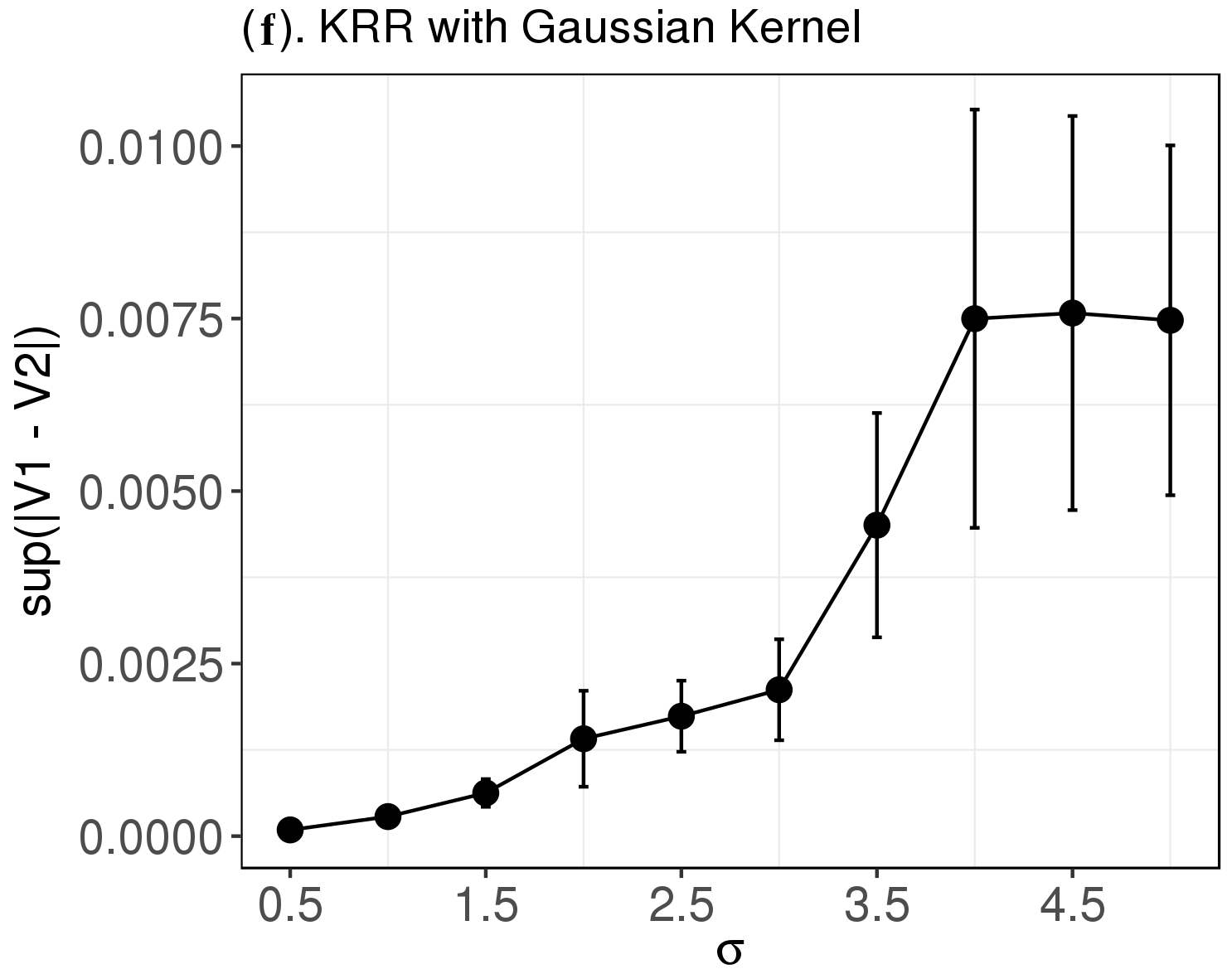}\\
\end{tabular}
  \caption{$(a)$, $(c)$, $(e)$: $\sup_{x} |V_1(x)-V_2(x)|$ for the polynomial kernel; $(b)$, $(d)$, $(f)$: $\sup_{x} |V_1(x)-V_2(x)|$ for the Gaussian kernel.  Error bars at $95\%$ confidence level for 200 repetitions of the experiments.}
\end{figure}

\section{Concluding Remarks}
There are several ways of extending this research. While we focused on kernel ridge regression, it would be interesting to propose a statistically and computationally efficient variance estimator for Gaussian processes as well. Additionally, 
currently in Assumption A1, we only considered sub-Gaussian random matrices, for theoretical convenience. However, the property in Assumption A1 might also hold for subsampled Fourier and Hadamard random matrices, but with a different relationship between $m$ and $s_\lambda$. For sub-Gaussian random matrices, we only need $m > s_\lambda$. But for Hadamard random matrices, $m > s_\lambda \log n $ is needed for estimation of the (mean) prediction error as in \cite{yang2015randomized}. This might also likely happen in our predictive variance. But note that our definition of $s_\lambda$ is different from \cite{yang2015randomized}. The analysis of different random matrices is appealing for future work. However, our general results on sub-Gaussian matrices should be seen as a necessary first step towards this endeavor.

\bibliographystyle{plainnat}
\bibliography{mml}

\section*{APPENDIX}\label{sec:app}
\subsection{Illustrative Application: Randomly Sketched Active Learning Algorithm}\label{sec:app}

Active learning has been successfully applied to classification as well as regression problems \cite{brinker2003incorporating}. Most active learning algorithms need to iteratively compute a score for each unlabeled samples. Specifically, the kernel ridge regression approach needs to evaluate the prediction variance for the unlabeled samples. As the size of the training data increases, the cost of computation increases cubically. An additional aspect that increases the computational cost is the use of cross validation to select the tuning parameters at each iteration, followed by computing the score for each unlabeled subject. 
Our randomly sketched active learning is aimed to reduce the computational cost for both model fitting and score calculation.

The main computational cost for randomly sketched KRR lies in computing the matrix multiplication of the sketch matrix and the kernel matrix. Suppose the current training set has $n_0$ data points and the projection dimension of the random matrix is $m$. The computational complexity of the matrix multiplication is of the order of $O(mn_0^2)$. In the next iteration, $n_s$ data points are added to the training set. Instead of calculating the matrix multiplication for all $n_0+n_s$ data points, we only need to calculate the entries corresponding to the updated data points. We partition the new kernel matrix as
\begin{equation*}
K = \begin{bmatrix}
K_1& K_{12}\\
K_{21} & K_2 
\end{bmatrix}
\end{equation*}
where $K_1 \in \mathbb{R}^{n_0 \times n_0}$ is the kernel matrix of the current training set, while $K_{12} \in \mathbb{R}^{n_0\times n_s}$, $K_{21} \in \mathbb{R}^{n_s\times n_0}$ and $K_2 \in \mathbb{R}^{n_s\times n_s}$. Correspondingly, we partition the new random projection matrix as,
\begin{equation*}
S = \begin{bmatrix}
S_1& S_{12}\\
S_{21} & S_2 
\end{bmatrix}
\end{equation*}
where $S_1 \in \mathbb{R}^{m_1 \times n_0}$ is the sketch matrix from the current step, while $S_{12} \in \mathbb{R}^{m_1 \times n_s}$, $S_{21} \in \mathbb{R}^{(m_2-m_1)\times n_0}$, $S_2 \in \mathbb{R}^{(m_2-m_1)\times n_s}$, and $m_1$ and $m_2$ are the projections dimension for the current and new sketch matrices correspondingly. Then the matrix multiplication can be written as,
\begin{equation}\label{eq:part}
\begin{aligned}
SK &= \begin{bmatrix}
K_1& K_{12}\\
K_{21} & K_2 
\end{bmatrix}\begin{bmatrix}
S_1& S_{12}\\
S_{21} & S_2 
\end{bmatrix} \\
& = \begin{bmatrix}
S_1 K_1 + S_{12}K_{21} & S_1 K_{12}+S_{12}K_{22} \\ 
S_{21}K_1+S_2K_{21} & S_{21}K_{12} + S_2 K_2
\end{bmatrix}
\end{aligned}
\end{equation} 
Since $S_1 K_1$ has already been calculated in the previous step, we only need to calculate the remaining terms. The computational complexity is thus reduced from $O(m(n_0+n_s)^2)$ to $O(mn_0 n_s)$. The size $n_0$ increases at each iteration, but the step size $n_s$ is fixed. Thus, $O(m n_0 n_s)$ is at most $O(mn_0)$. The reduction of computational complexity is significant for large training sets.
\begin{algorithm}[h]
   \caption{Active learning algorithm}
   \label{alg:example}
\begin{algorithmic}
   \STATE {\bfseries Input:} \\
   Initial training data set $S$\\
   Unlabeled data $U$
   \REPEAT
   \STATE \textbf{Step 1:} Calculate the projected kernel matrix using eq.(\ref{eq:part}).
   \STATE \textbf{Step 2:} Apply the randomly sketched kernel ridge regression to the training data $S$.
   \STATE \textbf{Step 3:} Calculate the randomly sketched prediction variance $V_2(x)$ for the samples in $U$ as in eq.(3.4).
   \STATE \textbf{Step 4:} Sample $n_s$ points based on the weight $V_2(x)$ and obtain the labels associated to them.
   \STATE \textbf{Step 5:} Add the sampled points to the training data set $S$ and remove them from unlabeled set $U$. 

   \UNTIL{The predefined convergence condition is satisfied}
   \STATE {\bfseries Output:}\\
   Final training set $S$.
\end{algorithmic}
\end{algorithm}

The difference between Algorithm 1 and the classical active learning algorithm is that we use $V_2(x)$ as weights to randomly sample data points from the unlabeled training data instead of deterministically selecting the data points with largest scores. If there is a small cluster of data with large scores, the deterministic method tends to add all of them into the training set initially. Suppose these data points are clustered together and outside of the majority of data points. Once they are all selected in the first few iterations, they may become the majority in the training set, and since the total size of labeled data is very small in early iterations, this will add an extra bias in the prediction. Thus we use the weighted random sampling strategy to ensure a substantial probability to select data points with large score while avoiding to select too many of them at once.  

\subsection{Experiments}\label{sec:append:sim}
In this section, we evaluate the performance of our proposed random projection approach. We run experiments on synthetic data as well as on real-world data sets.

\subsubsection{Confirming our theoretical contribution}
 Through synthetic experiments, we first verify the validity of our theoretical contribution (Theorem \ref{thm:main}). Here $500$ training samples were generated based on eq.(\ref{model}) with  $X \sim 1/2N(0.5,0.5) + 1/2 N(5,5)$. We use the Gaussian kernel function $K(x,x')=\exp{-\frac{(x-x')^2}{2\sigma^2}}$, where $\sigma=1$. Next, we generated $50$ testing samples following the same distribution. Here we generated a random projection matrix with Gaussian distributed entries, and the projection dimension is chosen as $m=c\sqrt{\log(n)}$, with $c=8,10,12$ respectively ($m=20,25,30$ approximately). We observe that, with the increase of $m$, the randomly projected variance performs similar to the original conditional variance, which confirmed the validity of our approach. Also, as for the computational time, the time for calculating the original variance for a new sample $x$ takes $4.654s$, but our proposed new randomly projected variance only takes $0.251s$, showing the practical advantage of our method.

\begin{figure}[h]\label{fig:compare}
  \centering
   \begin{tabular}{cc}
      \includegraphics[scale=0.34]{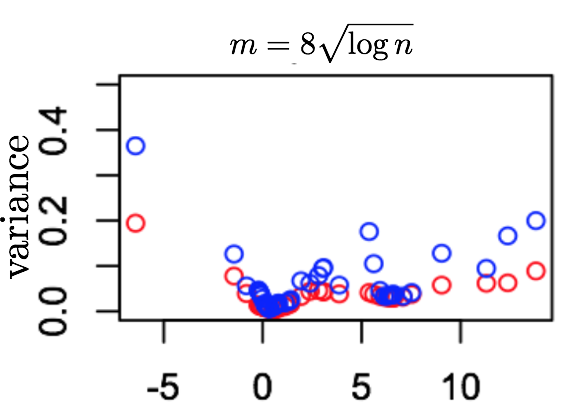}&
    \includegraphics[scale=0.34]{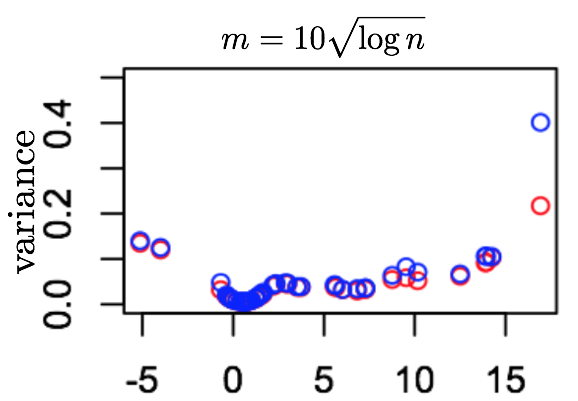}\\
    \includegraphics[scale=0.34]{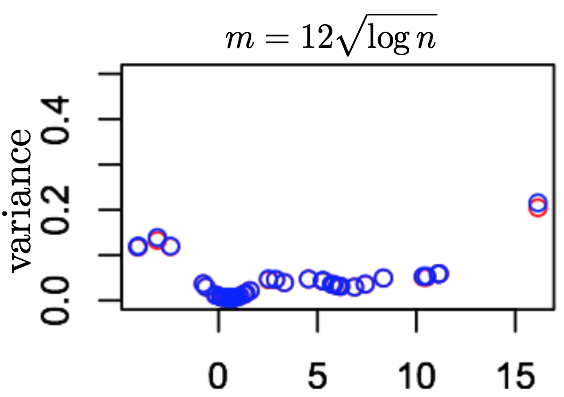}&\\
\end{tabular}
  \caption{Red dots represent the original variance for the new sample $x$, blue dots represent the randomly projected variance.}
\end{figure}

\subsubsection{Illustrative application on active learning: synthetic experiments}

Next, we illustrate the use of our variance estimator in active learning with synthetic data. (Appendix \ref{sec:app} provides details of a simple algorithm that uses of our variance estimator and attains $O(mn)$ time.) For comparison, we simulated $5000$ data points as the training set and $1000$ data points as the testing set. The initial training set was selected by randomly sampling $100$ data points from the training set. In the simulation settings, we use the Gaussian kernel,
\begin{equation*}
\mathcal{K}_{gau}(u,v) = e^{-\frac{1}{2h^2}(u-v)^2}
\end{equation*}
with bandwidth $h=0.25$. We report the mean squared error(MSE) at each iteration.

\textbf{Simulation Setting 1.} 
We simulate the predictor $X$ from a uniform distribution on $[0,1]$ and $f^{\ast}(x) = -1 + 2 x^2$. The response $y_i$ was generated as
$y_i = f^{\ast}(x_i) + \epsilon_i\;(i=1,\dots,n), $
where $\epsilon_i$ are i.i.d. standard Gaussian noise.

Gaussian random projection matrix is used in this setting, and we choose the sketch dimension $m=\lceil\log(n)\rceil$. As shown in Figure 2, the randomly sketched active learning algorithm has the smallest mean squared error after $30$ iteration. Also, the mean squared error of randomly sketched active learning algorithm converges as fast as the active learning with the original KRR and random sampling with original KRR.

\textbf{Simulation Setting 2.} 
We simulate the predictor $X$ from the following distribution
\begin{equation*}
x_i = \left\{
            \begin{array}{ll}
              \mbox{Unif}[0,1/2]\quad \mbox{if } i=1,\dots,k \\
              1 + z_i \quad \mbox{if } i=k+1,\dots,n 
            \end{array}
              \right.
\end{equation*}
where $z_i\sim N(0,1/n)$ and $k=\lceil \sqrt{n} \rceil$. For this experiment, we make  
$f^{\ast}(x) = -1 + 2 x^2$. The response $y_i$ was generated as
$y_i = f^{\ast}(x_i) + \epsilon_i \; (i=1,\dots,n) $, where $\epsilon_i$ are i.i.d. standard Gaussian noise.

\begin{figure}[h]\label{fig:sim}
  \centering
    \begin{tabular}{cc}
      \includegraphics[scale=0.24]{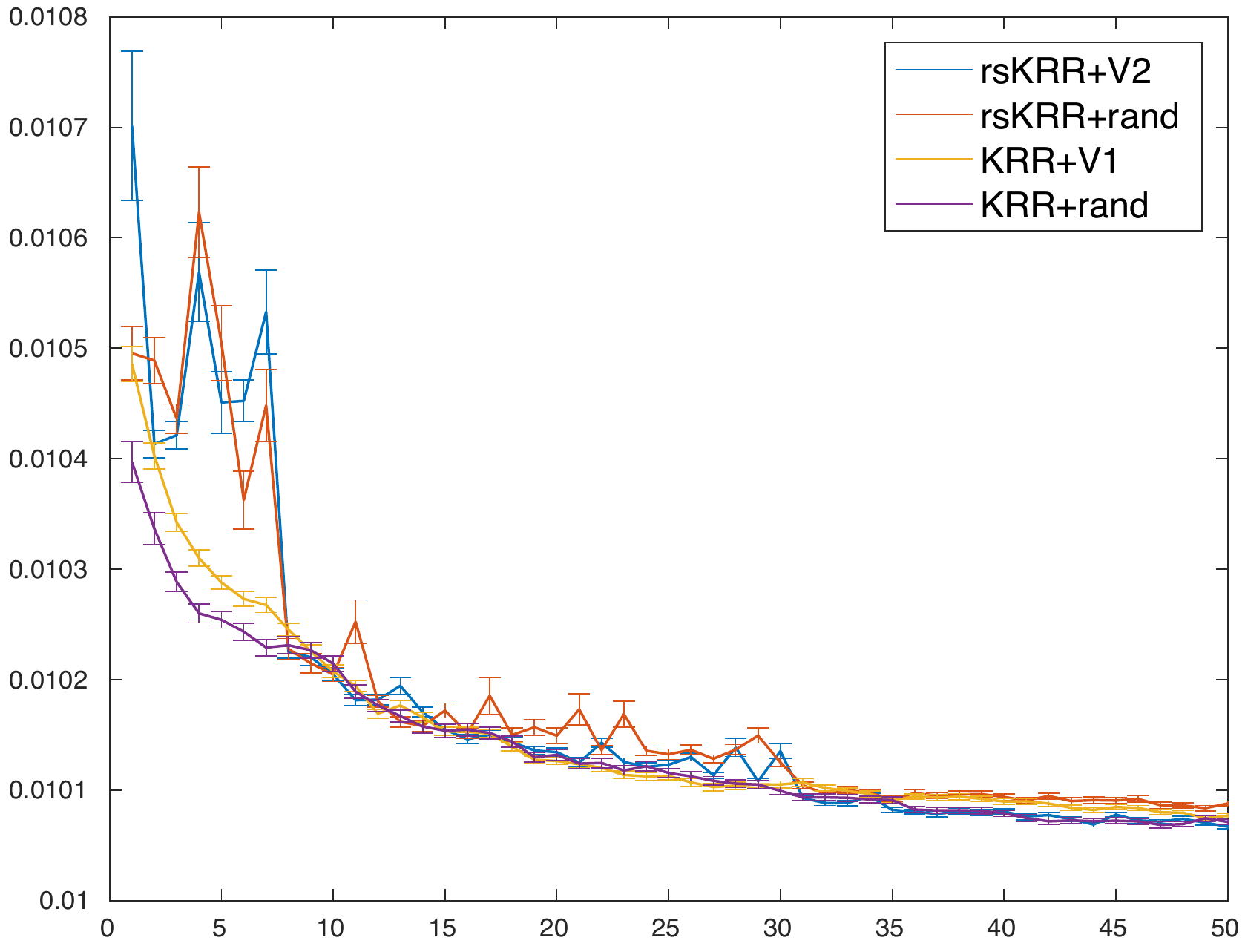}&
    \includegraphics[scale=0.24]{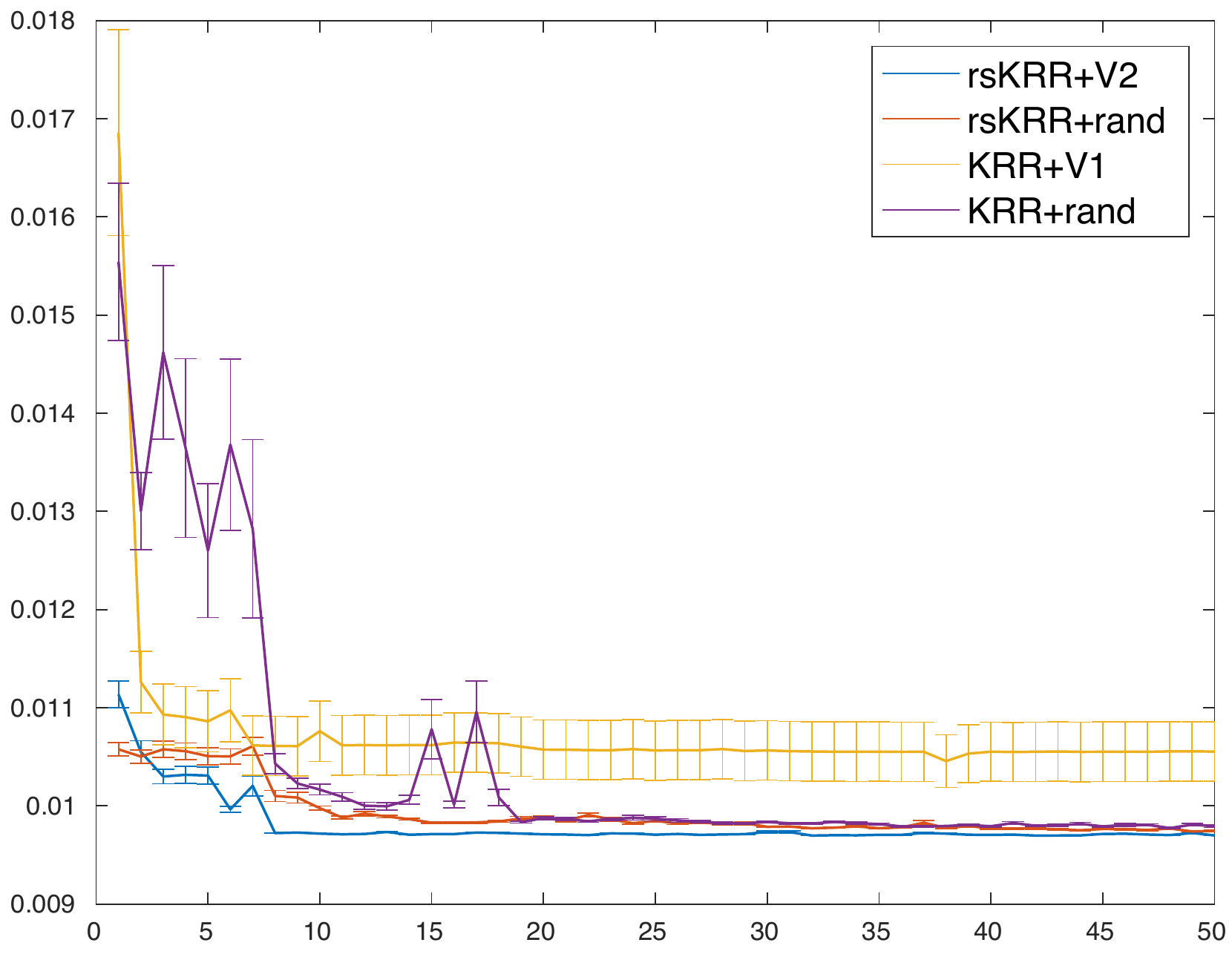}\\
    ~~~~~Simulation setting 1 &~~~~~ Simulation setting 2\\
\end{tabular}
  \caption{We compare active learning strategies with random sampling strategies under the original KRR, and the randomly sketched KRR respectively. Y-axes is the mean squared error. At each iteration, we add $30$ data points and show the predicted MSE of the four strategies in different colors:``rsKRR+V2'' denotes our randomly sketched active learning algorithm, ``KRR+V1'' denotes active learning with original KRR, ``KRR+rand'' denotes uniform random sampling with original KRR and ``rsKRR+rand'' denotes uniform random sampling with randomly sketched KRR. (Error bars at $95\%$ confidence level for $30$ repetitions of the experiments.)}
\end{figure}

Same as the Setting 1, we also use a Gaussian random matrix with sketch dimension $m=\lceil \log(n) \rceil$. As shown in Figure 3, the original active learning method shows faster convergence rate and achieves lower MSE after $50$ iterations compared to the random sampling algorithm. For this unevenly distributed data, it is unlikely to select the data outside the majority for the random sampling strategy. However, the minority data with large prediction variance tends to be selected by the active learning algorithm. Thus the randomly sketched active learning algorithm is comparable with active learning with the original KRR after $30$ iterations and converges to a similar MSE.

\subsubsection{Illustrative application on active learning: real-word experiments}
Next, we illustrate the use of our variance estimator in active learning with real-world data. (Appendix \ref{sec:app} provides details of a simple algorithm that uses of our variance estimator and attains $O(mn)$ time.)

\textbf{Flight Delay Data.} Here, we evaluate our randomly sketched active learning algorithm on the US flight dataset \cite{hensman2013gaussian} that contains up to 2 million points. We use a subset of the data with flight arrival and departure times for commercial flights in 2008. The flight delay was used as our response variable and we included 8 of the many variables from this dataset: the age of the aircraft, distance that needs to be covered, airtime, departure time, arrival time, day of the week, day of the month and month.

\begin{figure}[h]\label{fig:real}
\centering
   \begin{tabular}{cc}
      \includegraphics[scale=0.24]{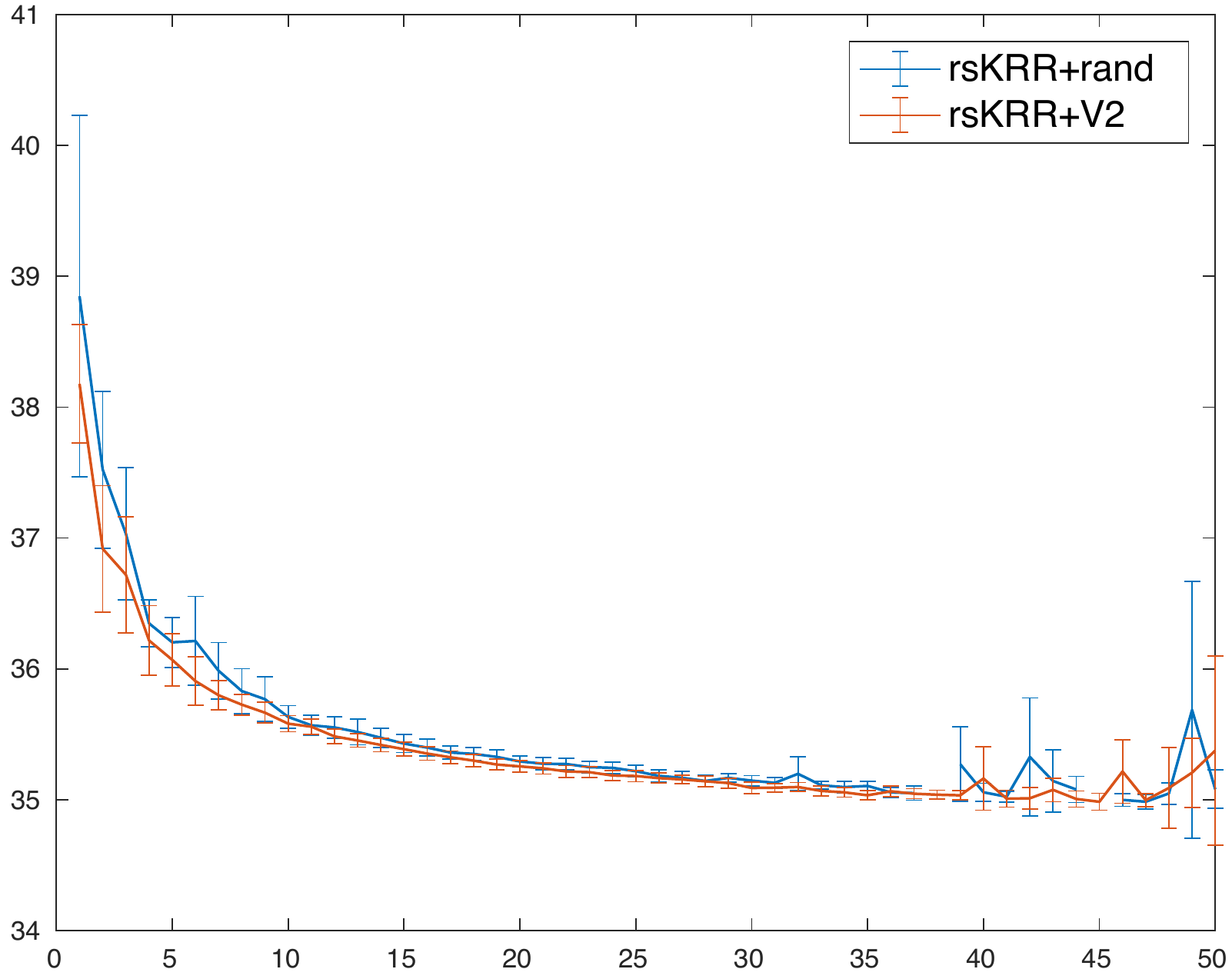}&
    \includegraphics[scale=0.24]{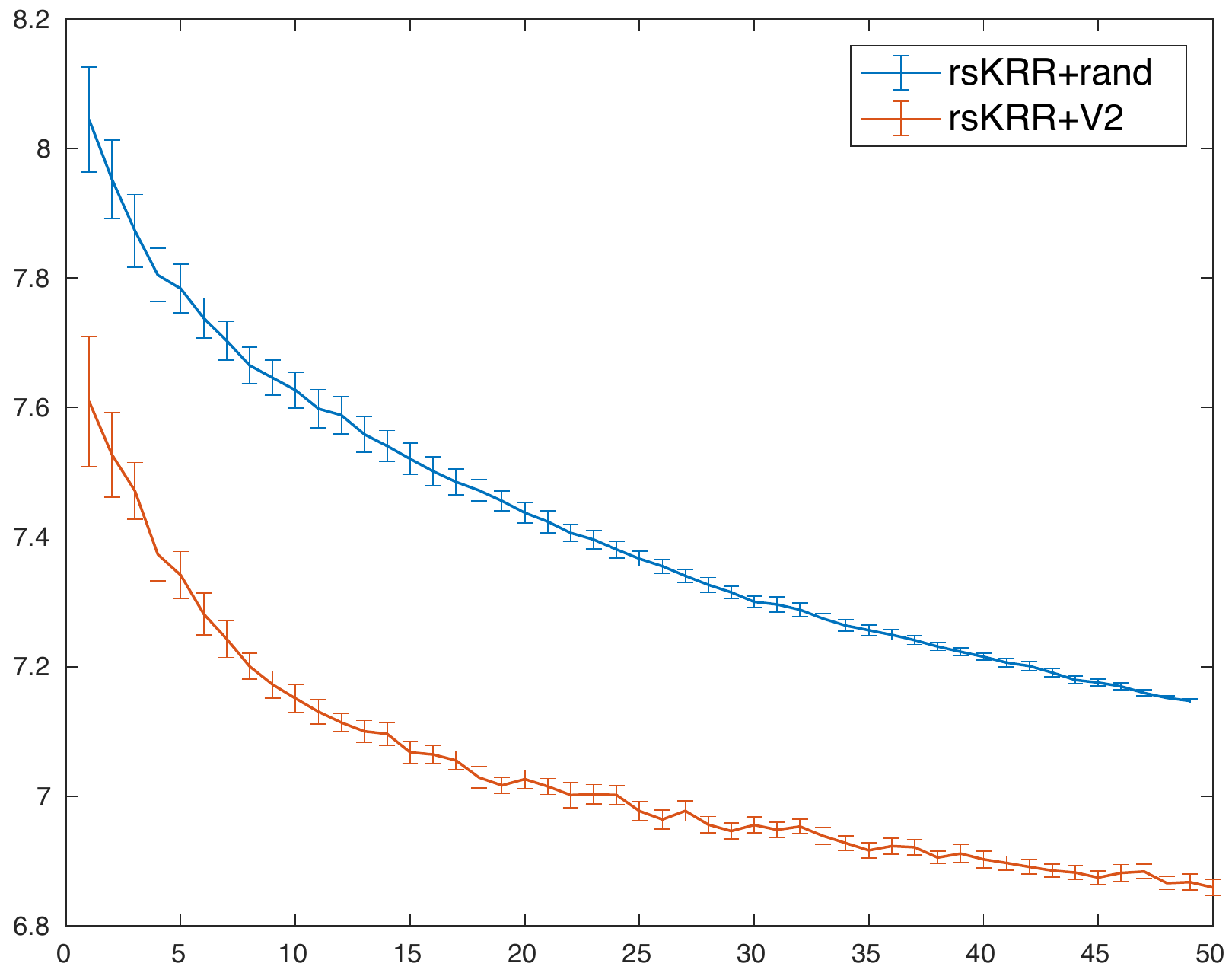}\\
    ~~~~~Flight delay data &~~~~~ World weather data\\
\end{tabular}
\caption{We compare active learning strategies with random sampling strategies using randomly sketched KRR. Y-axes is the root mean squared error. At each iteration, we add $1000$ data points and show the predicted rMSE of two strategies in different colors: ``rsKRR+V2'' denotes our randomly sketched active learning algorithm, ``rsKRR+rand'' denotes uniform random sampling with randomly sketched KRR. (Error bars at $95\%$ confidence level for 30 repetitions of the experiments.)}
\end{figure}

We randomly selected $60,000$ data points, using $50,000$ as the training set and $10,000$ as the testing set. We first randomly selected $1000$ data points as labeled data. Then we sequentially added $1000$ data points from the unlabeled training data at each iteration. We use the Gaussian random matrix with projection dimension $m=\lceil \log(n)\rceil$.  Here we only use the randomly sketched KRR since the computational cost and required RAM of the original KRR is too large. To compare the performance of active learning and uniform sampling, we calculate the RMSE(root mean squared error) 30 times using the prediction on the testing set. In Figure 4, the active learning algorithm achieves the RMSE of the full data faster than the uniform random sampling method. 

\textbf{World Weather Data.} In what follows, we examined our method on another real world dataset. The world weather dataset contains monthly measurements of temperature, precipitation, vapor, cloud cover, wet days and frost days from Jan 1990 to Dec 2002 on a $5\times 5$ degree grid that covers the entire world. In our experiments, the response variable is temperature. We use the Gaussian random matrix with projection dimension $m=\lceil \log(n)\rceil$. We use $10,000$ samples for training and $10,000$ samples for testing. We start with an initial set of $200$ labeled points, and add $200$ points at each iteration. As we can observe in Figure 5, our method compares favorably.

\end{document}